\documentclass[12pt,reqno,a4wide]{amsart}

\usepackage{tikz} 
\usepackage{calrsfs}
\usepackage{mathrsfs}
\usetikzlibrary{decorations.pathmorphing}
\usetikzlibrary{arrows.meta}
\usepgflibrary {shadings}
\usepackage{array}
\usepackage{xcolor}
\usepackage{comment}
\usetikzlibrary{decorations.pathreplacing,decorations.markings,snakes}

\usepackage{hyperref}

\usetikzlibrary{decorations.pathreplacing}
\usetikzlibrary{fit}

\usepackage{pgfplots}

\allowdisplaybreaks

\begin{document}
\bibliographystyle{plain}

%
%

	\title[Cornerless, peakless, valleyless Motzkin paths]
	{Cornerless, peakless, valleyless Motzkin paths (regular and skew) and applications to bar graphs}

	\author[H. Prodinger ]{Helmut Prodinger }
	\address{Department of Mathematics, University of Stellenbosch 7602, Stellenbosch, South Africa
	and
NITheCS (National Institute for
Theoretical and Computational Sciences), South Africa.}
	\email{hproding@sun.ac.za}


	\begin{abstract}
				
		Motzkin paths consist of up-steps, down-steps, horizontal steps, never go below the $x$-axis and return to the $x$-axis.
		Versions where the return to the $x$-axis isn't required are also considered.
		A path is peakless (valleyless) if $UD$ (if $DU$) never occurs. If it is both peakless and valleyless, it is called cornerless.
Deutsch and		Elizalde have linked cornerless Motzkin paths and bargraphs bijectly. Thus, instead of prefixes of bargraphs
one might consider prefixes of cornerless Motzkin paths. In this paper, this is extended by counting the occurrences of $UD$ resp.\  $DU$.

The concepts are extended to so-called skew Motzkin paths. Methods are generating functions and the kernel method to compute
explicit forms.

	\end{abstract}
	

\maketitle


\section{Introduction}

Motzkin  paths consist of steps $U=(1,1)$, $D=(1,-1)$, $H=(1,0)$, start at the origin, never go below the
$x$-axis and return to the $x$-axis. Motzkin meanders (partial Motzkin paths, prefixes of Motzkin paths) are similar, but must not
necessarily return to the $x$-axis. The letter $H$ stands for `horizontal', sometimes one finds the symbol
$F$ (`flat'). Motzkin meanders can be seen as prefixes of Motzkin paths, \cite{Blecheraeq}.
The notations `excursion' and `meander' were Philippe Flajolet's favorites, see e. g. \cite{BF}.

Any sequence labelled A*** in this paper refers to the encyclopedia of integer sequences \cite{OEIS}.
\begin{figure}[h]
	\begin{center}
		\begin{tikzpicture}[scale=0.4 ]
			
			\draw[thin,-latex] (0,0) -- (24,0);
			\draw[thin,-latex] (0,0) -- (0,4);		
			\draw [thick](0,0)-- (1,0) -- (2,1) --(3,1) --(4,1) --(5,2) --(6,2) --(7,2) --
			(8,1) --(9,0) --(10,1) --(11,2) --(12,3) --(13,3) --(14,2) -- (15,2) --(16,2) --(17,2) --(18,1) --
			(19,1) --(20,0) --(21,0) --(22,0) ;
			\node at (0,0) {$\bullet$};
			\node at (1,0) {$\bullet$};
			\node at (2,1) {$\bullet$};
			\node at (3,1) {$\bullet$};
			\node at (4,1) {$\bullet$};
			\node at (5,2) {$\bullet$};
			\node at (6,2) {$\bullet$};
			\node at (7,2) {$\bullet$};
			\node at (8,1) {$\bullet$};
			\node at (9,0) {$\bullet$};
			\node at (10,1) {$\bullet$};
			\node at (11,2) {$\bullet$};
			\node at (12,3) {$\bullet$};
			\node at (13,3) {$\bullet$};
			\node at (14,2) {$\bullet$};
			\node at (15,2) {$\bullet$};
			\node at (16,2) {$\bullet$};
			\node at (17,2) {$\bullet$};
			\node at (18,1) {$\bullet$};
			\node at (19,1) {$\bullet$};
			\node at (20,0) {$\bullet$};
			\node at (21,0) {$\bullet$};
			\node at (22,0) {$\bullet$};
		\end{tikzpicture}%
		\begin{tikzpicture}[scale=0.4 ]
			\draw[thin,-latex] (0,0) -- (15,0);
			\draw[thin,-latex] (0,0) -- (0,4);		
			\draw [thick](0,0)-- (1,0) -- (2,1) --(3,1) --(4,1) --(5,2) --(6,2) --(7,2) --
			(8,1) --(9,0) --(10,1) --(11,2) --(12,3) --(13,3) --(14,2);
			\node at (0,0) {$\bullet$};
			\node at (1,0) {$\bullet$};
			\node at (2,1) {$\bullet$};
			\node at (3,1) {$\bullet$};
			\node at (4,1) {$\bullet$};
			\node at (5,2) {$\bullet$};
			\node at (6,2) {$\bullet$};
			\node at (7,2) {$\bullet$};
			\node at (8,1) {$\bullet$};
			\node at (9,0) {$\bullet$};
			\node at (10,1) {$\bullet$};
			\node at (11,2) {$\bullet$};
			\node at (12,3) {$\bullet$};
			\node at (13,3) {$\bullet$};
			\node at (14,2) {$\bullet$};
		\end{tikzpicture}
		
	\end{center}
	\caption{A Motzkin excursion and a Motzkin meander (ending at level $2$).}
	\label{motz1}
\end{figure}
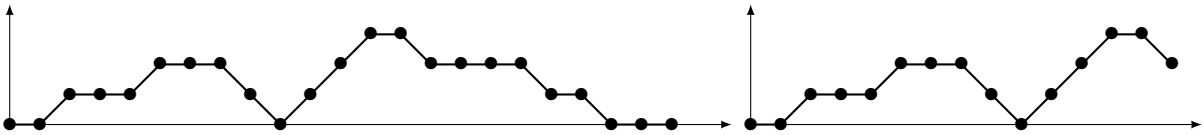

Motzkin paths can be recognized by a graph (automaton); if one ends at level 0, it is a Motzkin path. If it ends at level $j$ or if $j\ge0$ can take any value, we speak 
about prefixes of Motzkin paths.
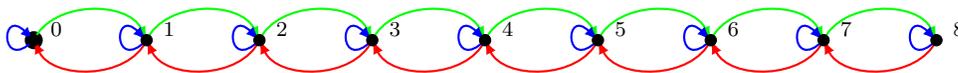
\begin{figure}[h]

	\begin{center}
		\begin{tikzpicture}[scale=1.5,main node/.style={circle,draw,font=\Large\bfseries}]

			\fill (0,-1) circle (0.08cm);

			\foreach \x in {0,1,2,3,4,5,6,7}
			{
				\draw[thick, -latex,green] (\x,-1) to  [out=65,in=115](\x+1,-1);

				\draw[ thick,red, -latex] (\x+1,-1) to  [out=-115,in=-65](\x,-1);    
				\draw[thick,latex-,blue ] (\x,-1)  ..  controls (\x-0.3,-1+0.3) and  (\x-0.3,-1-0.3) .. (\x,-1) ;	
				
				\node at  (\x+0.2,0.1-1){\tiny$\x$};
			}  
			
			\foreach \x in {8}
			{
				\draw[thick,latex-,blue ] (\x,-1)  ..  controls (\x-0.3,-1+0.3) and  (\x-0.3,-1-0.3) .. (\x,-1) ;
			}

			\node at  (8+0.2,0.1-1){\tiny$8$};
			
			\foreach \x in {0,1,2,3,4,5,6,7,8}
			{
				\draw (\x,-1) circle (0.05cm);
				\fill (\x,-1) circle (0.05cm);
			}
		\end{tikzpicture}
	\end{center}
	\caption{Graph (automaton) to recognize  Motzkin paths  }
	\label{purpelUDDU}
\end{figure}

If an up-step and a down-step are added at the beginning resp.\ end, we get \emph{elevated} Motzkin paths. The use of them will become clear later when we discuss the connections to bargraphs. They can also be recognized by an automaton with a special state.

\begin{figure}[h]

	\begin{center}
		\begin{tikzpicture}[scale=1.5,main node/.style={circle,draw,font=\Large\bfseries}]

			\fill (-1.5,-1) circle (0.08cm);
			\draw[ultra thick, -latex,teal] (-1.5,-1) to  [out=65,in=115](0,-1);
			\draw[ultra thick,purple, -latex] (0,-1) to [in=0,out=-120](-1.5,-1.5);    
			
						\fill[gray] (-1.5,-1.5) circle (0.08cm);
			\foreach \x in {0,1,2,3,4,5,6}
			{
				\draw[thick, -latex,green] (\x,-1) to  [out=65,in=115](\x+1,-1);

				\draw[ thick,red, -latex] (\x+1,-1) to  [out=-115,in=-65](\x,-1);    
				\draw[thick,latex-,blue ] (\x,-1)  ..  controls (\x-0.3,-1+0.3) and  (\x-0.3,-1-0.3) .. (\x,-1) ;	
				
				
			}  
			
			\node at  (-1.7,0.2-1){\scriptsize{start}};
						\node at  (-1.65,0.2-1.53){\scriptsize{end}};
			
			\foreach \x in {0,1,2,3,4,5,6,7}
			{	
				\node at  (\x+0.2,0.1-1) {\tiny\x};
			}  
			
			\foreach \x in {7}
			{
				\draw[thick,latex-,blue ] (\x,-1)  ..  controls (\x-0.3,-1+0.3) and  (\x-0.3,-1-0.3) .. (\x,-1) ;
			}

			
			\foreach \x in {0,1,2,3,4,5,6,7}
			{
				\draw (\x,-1) circle (0.05cm);
				\fill (\x,-1) circle (0.05cm);
			}
		\end{tikzpicture}
	\end{center}
	\caption{Graph (automaton) to recognize elevated Motzkin paths  }
	\label{purpelelevate}
\end{figure}
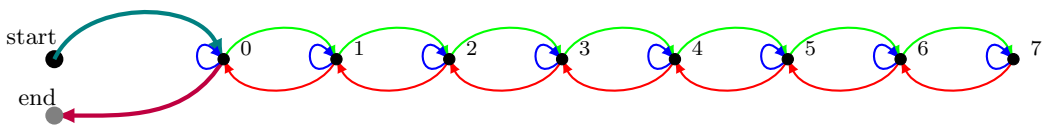

Motzkin paths without horizontal steps, thus only comprising up-steps and down-steps are called \emph{Dyck} paths. This family is much older and more prominent; their enumeration leads to \emph{Catalan} numbers which are the subject of a whole book, see \cite{Stanley-book}.

A variation of Dyck paths was introduced in \cite{Deutsch-italy}; a third step $L=(-1,-1)$ is introduced, but the resulting path is not allowed to intersect itself, which
can be described by forbidding the contiguous subwords $UL$ and $LU$. This variation is called \emph{skew} Dyck paths. 

The concept can be lifted to the instance of Motzkin paths, leading to \emph{skew} Motzkin paths, again forbidding the contiguous subwords $UL$ and $LU$. 
I have seen this concept first in \cite{Qing}. Some progress has been reported in \cite{garden}.

\begin{figure}[h]

	\begin{center}
		\begin{tikzpicture}[scale=1.6]
			\draw (0,0) circle (0.08cm);
			\fill (0,0) circle (0.08cm);
			
			\foreach \x in {0,1,2,3,4,5,6,7,8}
			{
				\draw (\x,0) circle (0.05cm);
				\fill (\x,0) circle (0.05cm);
			}
			
			\foreach \x in {0,1,2,3,4,5,6,7,8}
			{
				\draw (\x,-1) circle (0.05cm);
				\fill (\x,-1) circle (0.05cm);
			}
			
			\foreach \x in {0,1,2,3,4,5,6,7,8}
			{
				\draw (\x,-2) circle (0.05cm);
				\fill[cyan] (\x,-2) circle (0.05cm);
			}
			\foreach \x in {0,1,2,3,4,5,6,7,8}
			{
				\draw (\x,-3) circle (0.05cm);
				\fill[red] (\x,-3) circle (0.05cm);
			}

			\foreach \x in {0,1,2,3,4,5,6,7}
			{
				\draw[ thick,-latex] (\x,0) -- (\x+1,0);
				
			}
			
			\foreach \x in {0,1,2,3,4,5,6,7}
			{
				\draw[ thick,-latex] (\x+1,0) to [out=200, in =70] (\x,-1);
				\draw[ thick,-latex] (\x,-1) to [out=30, in =-110] (\x+1,0);
				\draw[ thick,-latex] (\x+1,-1) to (\x,-1);
			}
			
			\foreach \x in {0,1,2,3,4,5,6,7,8}
			{
				\draw[ thick,cyan,-latex] (\x,0) to  [out=-55, in =55] (\x,-2);
				\draw[ thick,cyan,-latex] (\x,-1) to   (\x,-2);
				\draw[ thick,cyan,-latex] (\x,-3) to   (\x,-2);
			}

			\foreach \x in {0,1,2,3,4,5,6,7}
			{
				\draw[ thick,red,-latex] (\x+1,-1) to  (\x,-3);
				\draw[ thick,red,-latex] (\x+1,-2) to[out=-110, in =30]   (\x,-3);
				\draw[ thick,red,-latex] (\x+1,-3) to   (\x,-3);
			}
			
			\foreach \x in {0,1,2,3,4,5,6,7}
			{
				\draw[ thick,-latex] (\x+1,-3) to    (\x,-1);
				
			}
			\foreach \x in {0,1,2,3,4,5,6,7,8}
			{
				\draw[ thick,cyan ] (\x,-2) to  [out=100, in =80]  (\x-0.2,-2);
				\draw[ thick,cyan,-latex ] (\x-0.2,-2 ) to  [out=-60, in =200] (\x,-2) ;	
			}
			
			\foreach \x in {0,1,2,3,4,5,6,7}
			{
				\draw[ thick,-latex] (\x,-2) to [out=30, in =-90]     (\x+1,-0);
				
			}
			
			\foreach \x in {0,1,2,3,4,5,6,7}
			{
				\draw[ thick,-latex] (\x+1,-2) to [out=100, in =-30]     (\x,-1);
				
			}
			
		\end{tikzpicture}
	\end{center}
	\caption{Four layers of states according to the type of steps leading to them. Traditional up-steps and down-steps are black, 
		level-steps are blue, and left steps are red.}
	\label{schoas}
\end{figure}
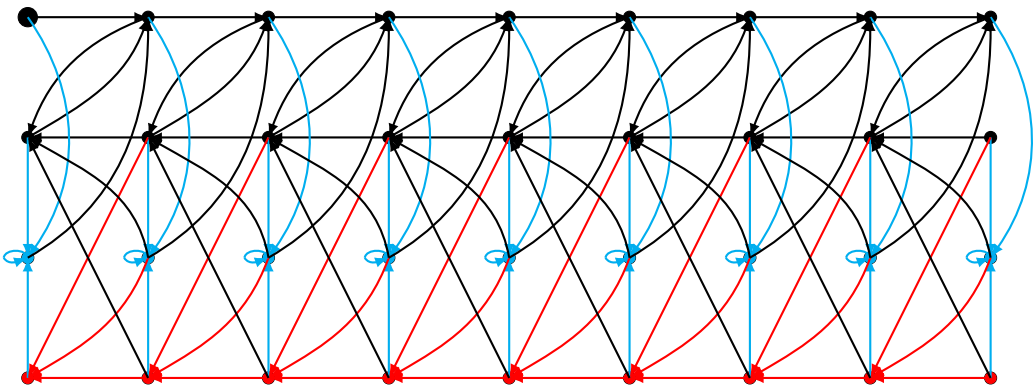

The situation is best described by a graph (state-diagram); see Figure~\ref{schoas}.

A Motzkin path (or a prefix of it) is called \emph{peakless} if a $U$-step is never followed by a $D$-step, or it is does not contain the contiguous subword $UD$.
A Motzkin path (or a prefix of it) is called \emph{valleyless} if a $U$-step is never followed by a $D$-step, or it is does not contain the contiguous subword $UD$.
A Motzkin path (or a prefix of it) is called \emph{cornerless} if it is peakless \emph{end} valleyless. Both contiguous subwords $DU$ and $UD$ are forbidden.

For the reader's convenience, here is an automaton to recognize peakless Motzkin paths (borrowed from \cite{prodinger-peakless}).

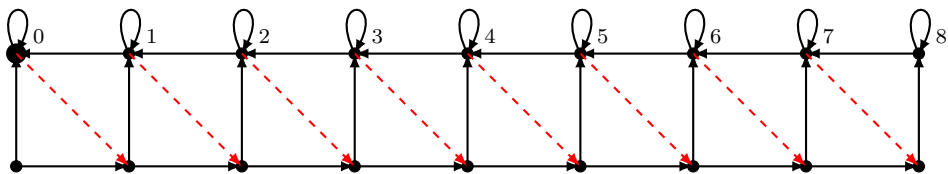
\begin{figure}[h]

	\begin{center}
		\begin{tikzpicture}[scale=1.5,main node/.style={circle,draw,font=\Large\bfseries}]

			\foreach \x in {0,1,2,3,4,5,6,7,8}
			{
				\draw (\x,0) circle (0.05cm);
				\fill (\x,0) circle (0.05cm);
				\draw (\x,-1) circle (0.05cm);
				\fill (\x,-1) circle (0.05cm);
			}

			\fill (0,0) circle (0.09cm);

			\foreach \x in {0,2,4,6}
			{
			}
			\foreach \x in {0,...,8}
			{
				\draw[thick,-latex ] (\x,0)  ..  controls (\x-0.25,0.5) and  (\x+0.25,0.5) .. (\x,0) ;	
			}

			\foreach \x in {0,1,2,3,4,5,6,7}
			{
				\draw[thick, latex-] (\x,0) to  (\x+1,0);	
				\draw[thick, -latex] (\x,-1) to  (\x+1,-1);	
				\draw[thick,red,dashed, -latex] (\x,0) to  (\x+1,-1);	
				\node at  (\x+0.2,0.15){\tiny$\x$};
			}			
			
			\foreach \x in {0,1,2,3,4,5,6,7,8}
			{
				\draw[thick, -latex] (\x,-1) to  (\x,0);	
			}

			\node at  (8+0.2,0.15){\tiny$8$};

		\end{tikzpicture}
	\end{center}
	\caption{Graph (automaton) to recognize peakless Motzkin paths; only the first few states are shown. Starting  at the origin and ending at nodes labelled 0 corresponds to Motzkin paths, and
		ending at a node labelled $k$ to  a path that ends at level $k$. }
	\label{purpel}
\end{figure}

And, again for convenience, here is an automaton to recognize cornerless Motzkin paths, borrowed from \cite{garden}.

\begin{figure}[h]

	\begin{center}
		\begin{tikzpicture}[scale=1.5,main node/.style={circle,draw,font=\Large\bfseries}]

			\fill (0,-1) circle (0.08cm);

			\foreach \x in {0,1,2,3,4,5,6,7}
			{
				\draw[thick, -latex,green] (\x,0-1) to  (\x+1,0);   
				\draw[ thick, latex-,green] (\x+1,0) to  (\x,0);

				\draw[ thick,red, -latex] (\x+1,-1) to  (\x,-2);    
				\draw[thick,latex-,blue ] (\x,-1)  ..  controls (\x-0.3,-1+0.3) and  (\x-0.3,-1-0.3) .. (\x,-1) ;	
				\draw[ thick,red, - latex] (\x+1,-2)   to  (\x,-2);    
				
				\node at  (\x+0.2,0.1){\tiny$\x$};
			}  
			
			\foreach \x in {8}
			{
				\draw[thick,latex-,blue ] (\x,-1)  ..  controls (\x-0.3,-1+0.3) and  (\x-0.3,-1-0.3) .. (\x,-1) ;
			}
			
			\foreach \x in {0,1,2,3,4,5,6,7,8}
			{
				
				\draw[ thick,blue, -latex] (\x,0)to  (\x,-1);    
				\draw[ thick,blue, -latex] (\x,-2)to  (\x,-1);    
				
			}                        
			
			\node at  (8+0.2,0.1){\tiny$8$};
			
			\foreach \x in {0,1,2,3,4,5,6,7,8}
			{
				\draw (\x,0) circle (0.05cm);
				\fill (\x,0) circle (0.05cm);
				\draw (\x,-1) circle (0.05cm);
				\fill (\x,-1) circle (0.05cm);
				\draw (\x,-2) circle (0.05cm);
				\fill (\x,-2) circle (0.05cm);
			}
		\end{tikzpicture}
	\end{center}
	\caption{Graph (automaton) to recognize  Motzkin paths with forbidden subwords UD and DU.}
	\label{purpelUDDU}
\end{figure}
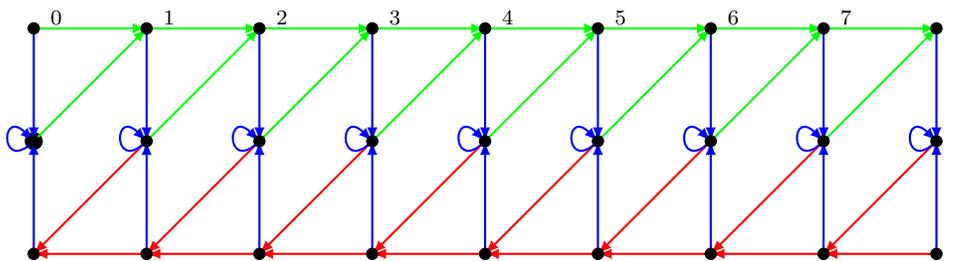

A main part of this paper will be the enumeration of Motzkin paths where peaks resp.\ valleys \emph{may} occur but are \emph{counted} using
extra variables $\sigma$ and $\tau$. If $\sigma=0$ we are in the peakless domain, if $\tau=0$, we  are in the valleyless domain, if both variables are zero, we
are in the cornerless world, and if $\sigma=\tau=1$, we have ordinary Motzkin paths. 
In a last section, the skew Motzkin paths are analyzed with respect to peaks resp.\ valleys. Recall that $UL$ and $LU$ are still not allowed.

\bigskip

We describe a major application of cornerless Motzkin paths, which deals with bargraphs;
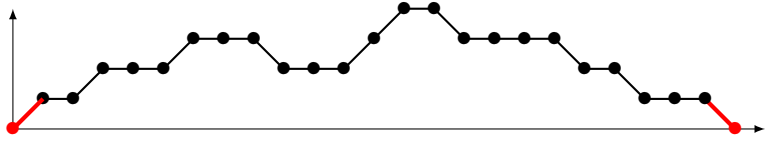
\begin{figure}[h]
	\begin{center}
		\begin{tikzpicture}[scale=0.4 ]
			
			\draw[thin,-latex] (-1,-1) -- (24,-1);
			\draw[thin,-latex] (-1,-1) -- (-1,3);		
			\draw [thick](0,0)-- (1,0) -- (2,1) --(3,1) --(4,1) --(5,2) --(6,2) --(7,2) --
			(8,1) --(9,1) --(10,1) --(11,2) --(12,3) --(13,3) --(14,2) -- (15,2) --(16,2) --(17,2) --(18,1) --
			(19,1) --(20,0) --(21,0) --(22,0) ;
			\draw (0,0) circle (0.08cm);
			\node at (0,0) {$\bullet$};
			\node[red] at (-1,-1) {$\bullet$};
						\node[red] at (23,-1) {$\bullet$};
						\draw [red, ultra thick](0,0)-- (-1,-1);
												\draw [red, ultra thick](22,0)-- (23,-1);
			\node at (1,0) {$\bullet$};
			\node at (2,1) {$\bullet$};
			\node at (3,1) {$\bullet$};
			\node at (4,1) {$\bullet$};
			\node at (5,2) {$\bullet$};
			\node at (6,2) {$\bullet$};
			\node at (7,2) {$\bullet$};
			\node at (8,1) {$\bullet$};
			\node at (9,1) {$\bullet$};
			\node at (10,1) {$\bullet$};
			\node at (11,2) {$\bullet$};
			\node at (12,3) {$\bullet$};
			\node at (13,3) {$\bullet$};
			\node at (14,2) {$\bullet$};
			\node at (15,2) {$\bullet$};
			\node at (16,2) {$\bullet$};
			\node at (17,2) {$\bullet$};
			\node at (18,1) {$\bullet$};
			\node at (19,1) {$\bullet$};
			\node at (20,0) {$\bullet$};
			\node at (21,0) {$\bullet$};
			\node at (22,0) {$\bullet$};
		\end{tikzpicture}%
		
	\end{center}
	\caption{Elevated cornerless Motzkin path}
	\label{motz1}
\end{figure}

Bargraphs are non-intersecting lattice path in $\mathbb{N}\times \mathbb{N}$
 with three allowable steps: up $U=(0, 1)$, down step $D=(0, -1)$ and horizontal step  $H=(1, 0)$. An up-step may
not immediately follow a down step nor vice versa; thus, $UD$ and $DU$ as contiguous subwords are forbidden.
 A bargraph starts at the
origin with an up step and terminates immediately upon return to the $x$-axis. Thus, the $x$-axis is only visited at the beginning and the end
of the bargraph, not inbetween.
We consider bargraphs defined by `semi-perimeter' (number of $U$-steps plus number of $H$-steps). See \cite{Blecheraeq} for more details;
Blecher and Knopfmacher analyze the progress of walking along a bargraph, but being allowed to stop anywhere. This results in \emph{prefixes} of bargraphs.

Instead of following such an approach directly, as Blecher and Knopfmacher did, we consider the equivalend concept of cornerless Motzkin paths. 
The bijection~\cite{DeutschElizalde} goes from cornerless Motzkin paths to elevated paths, by attaching a $U$-step in the beginning and a $D$-step at the end. Then we imagine the picture rotated by 45 degrees, and translating $U$ into $(0,1)$ leaving $H$ as $(1,0)$ and translating $D$ into $(0,-1)$. The pictures that are provided explain this easily.

For convenience, we ignore the additional steps that elevates the path, although this could be handled without too much effort. Recall that our own paper \cite{valerie} covers  several other instances where both, $UD$ and $DU$ are forbidden. The paper \cite{valerie} contains references to earlier efforts on the subject, notably by Valerie Roitner who included the material in her Ph.D. thesis~\cite{Vally}.

Here, we go one step further: we allow $UD$ but use a variable $\sigma$ to count often it occured, likewise we use $\tau$ to count the occurrences of $DU$. Eventually we are interested in the count $(n,j,k,\ell)$ where $n$ is the number of steps of the prefix of a Motzkin path, $j$ is the height of the end point, $k$ is the number of of $UD$'s and
$\ell$ is the number of $DU$'s. In generating functions, we use the variables $z,u,\sigma,\tau$. For convenience, we mostly write $\Phi(u)$  for what is really $\Phi(u;z,\sigma,\tau)$. ($\Phi$ is just a generic name of a generating function.)

\begin{figure}[h]
	\begin{center}
		\begin{tikzpicture}[scale=0.4 ]
			
			\draw[thin,-latex] (-1,-1) -- (14,-1);
			\draw[thin,-latex] (-1,-1) -- (-1,5);		
			\draw [thick](-1,0)-- (0,0) -- (0,1) --(1,1) --(2,1) --(2,2) --(4,2) --(4,1) --
			(4,1) --(6,1) --(6,3) --(7,3) --(7,2) --(10,2) --(10,1) -- (11,1) --(11,0) --(13,0);
			\draw (0,0) circle (0.08cm);
			\node at (-1,0) {$\bullet$};
			\node[red] at (-1,-1) {$\bullet$};
			\node[red] at (13,-1) {$\bullet$};
			\draw [red, ultra thick](-1,0)-- (-1,-1);
			\draw [red, ultra thick](13,0)-- (13,-1);
			\node at (0,0) {$\bullet$};
			\node at (0,1) {$\bullet$};
			\node at (1,1) {$\bullet$};
			\node at (2,1) {$\bullet$};
			\node at (2,2) {$\bullet$};
			\node at (3,2) {$\bullet$};
			\node at (4,2) {$\bullet$};
			\node at (4,1) {$\bullet$};
			\node at (5,1) {$\bullet$};
			\node at (6,1) {$\bullet$};
			\node at (6,3) {$\bullet$};
			\node at (6,2) {$\bullet$};
			\node at (7,3) {$\bullet$};
			\node at (7,2) {$\bullet$};
			\node at (8,2) {$\bullet$};
			\node at (9,2) {$\bullet$};
			\node at (10,2) {$\bullet$};
			\node at (10,1) {$\bullet$};
			\node at (11,1) {$\bullet$};
			\node at (11,0) {$\bullet$};
			\node at (12,0) {$\bullet$};
			\node at (13,0) {$\bullet$};
		\end{tikzpicture}\qquad%
		\begin{tikzpicture}[scale=0.4 ]
			
			\draw[thin,-latex] (-1,-1) -- (14,-1);
			\draw[thin,-latex] (-1,-1) -- (-1,4);		
			\draw [thick](-1,0)-- (0,0);
			 \draw [thick] (0,1) --(0,-1);
			  \draw [thick](1,-1) --(1,1);
			  \draw [thick](2,-1) --(2,2);
			   \draw [thick](0,1) --(2,1);
			   \draw [thick](3,-1) --(3,2);
			   			   \draw [thick](3,1) --(3,2);
			   			   			   \draw [thick](4,-1) --(4,2);
			   			   			   \draw [thick](4,2) --(2,2);
   			   			   			   			   \draw [thick](5,-1) --(5,1);
   			   			   			   			      			   			   			   			   \draw [thick](6,-1) --(6,3);
   			   			   			   			      			   			   			   			    \draw [thick](7,-1) --(7,3);
   			   			   			   			      			   			   			   			        \draw [thick](6,3) --(7,3);
\draw [thick](7,2) --(10,2);
\draw [thick](10,1) --(11,1);
\draw [thick](13,0) --(11,0);
 \draw [thick](8,-1) --(8,2);
  \draw [thick](9,-1) --(9,2);
   \draw [thick](10,-1) --(10,2);
    \draw [thick](11,-1) --(11,1);
    \draw [thick](12,-1) --(12,0);
    \draw [thick](13,-1) --(13,0);
    \draw [thick](4,1) --(6,1);
			\draw [  thick](-1,0)-- (-1,-1);
			\draw [  thick](13,0)-- (13,-1);
		\end{tikzpicture}%
	\end{center}
	\caption{After rotation; the bargraph}
	\label{motz1}
\end{figure}
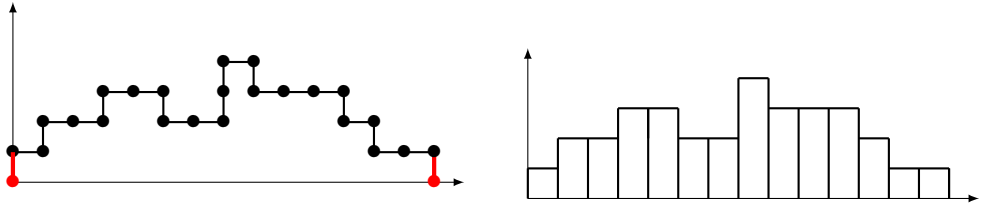
This important bijection is due to \cite{DeutschElizalde}.
	Here (Figure~\ref{motz3}) are some small examples:
\begin{figure}[h]
\label{motz3}
	\begin{center}
		\begin{tikzpicture}[scale=0.3 ]
						\draw[thin,-latex] (-1,-1) -- (7,-1);
			\draw[thin,-latex] (-1,-1) -- (-1,3);		
			\draw [thick] (0,0) -- (1,1) --(3,1) --(4,0);
			\draw (0,0) circle (0.08cm);
			\node at (0,0) {$\bullet$};
						\node at (1,1) {$\bullet$};
									\node at (2,1) {$\bullet$};
												\node at (3,1) {$\bullet$};
															\node at (4,0) {$\bullet$};
			\node[red] at (-1,-1) {$\bullet$};
			\node[red] at (5,-1) {$\bullet$};
			\draw [red, ultra thick](-1,-1)-- (0,-0);
			\draw [red, ultra thick](4,0)-- (5,-1);
		\end{tikzpicture}%
		\begin{tikzpicture}[scale=0.3,x=1cm ]
	\draw[thin,-latex] (-1,-1) -- (7,-1);
	\draw[thin,-latex] (-1,-1) -- (-1,3);		
	\draw [thick] (0,0) -- (1,0)-- (2,1) --(3,1) --(4,0);
	\draw (0,0) circle (0.08cm);
	\node at (0,0) {$\bullet$};
	\node at (1,0) {$\bullet$};
	\node at (2,1) {$\bullet$};
	\node at (3,1) {$\bullet$};
	\node at (4,0) {$\bullet$};
	\node[red] at (-1,-1) {$\bullet$};
	\node[red] at (5,-1) {$\bullet$};
	\draw [red, ultra thick](-1,-1)-- (0,-0);
	\draw [red, ultra thick](4,0)-- (5,-1);
\end{tikzpicture}%
		\begin{tikzpicture}[scale=0.3,x=1cm ]
	\draw[thin,-latex] (-1,-1) -- (7,-1);
	\draw[thin,-latex] (-1,-1) -- (-1,3);		
	\draw [thick] (0,0) -- (1,1) --(2,1)--(3,0) --(4,0);
	\draw (0,0) circle (0.08cm);
	\node at (0,0) {$\bullet$};
	\node at (1,1) {$\bullet$};
	\node at (2,1) {$\bullet$};
	\node at (3,0) {$\bullet$};
	\node at (4,0) {$\bullet$};
	\node[red] at (-1,-1) {$\bullet$};
	\node[red] at (5,-1) {$\bullet$};
	\draw [red, ultra thick](-1,-1)-- (0,-0);
	\draw [red, ultra thick](4,0)-- (5,-1);
\end{tikzpicture}%
		\begin{tikzpicture}[scale=0.3,x=1cm ]
	\draw[thin,-latex] (-1,-1) -- (7,-1);
	\draw[thin,-latex] (-1,-1) -- (-1,3);		
	\draw [thick] (0,0) -- (4,0);
	\draw (0,0) circle (0.08cm);
	\node at (0,0) {$\bullet$};
	\node at (1,0) {$\bullet$};
	\node at (2,0) {$\bullet$};
	\node at (3,0) {$\bullet$};
	\node at (4,0) {$\bullet$};
	\node[red] at (-1,-1) {$\bullet$};
	\node[red] at (5,-1) {$\bullet$};
	\draw [red, ultra thick](-1,-1)-- (0,-0);
	\draw [red, ultra thick](4,0)-- (5,-1);
\end{tikzpicture}

		\begin{tikzpicture}[scale=0.4,x=1cm ]
	\draw[thin,-latex] (-1,-1) -- (3,-1);
	\draw[thin,-latex] (-1,-1) -- (-1,2);		
	
	\draw[thick] (-1,-1) -- (-1,1);
	\draw[thick] (-1,1) -- (1,1);		
	\draw[thick] (1,1) -- (1,-1);		
	\draw[thick] (0,1) -- (0,-1);		
\end{tikzpicture}\qquad
		\begin{tikzpicture}[scale=0.4,x=1cm ]
	\draw[thin,-latex] (-1,-1) -- (3,-1);
	\draw[thin,-latex] (-1,-1) -- (-1,2);		
	
	\draw[thick] (-1,-1) -- (-1,0);
	\draw[thick] (-1,0) -- (0,0);		
	\draw[thick] (0,-1) -- (0,1);		
	\draw[thick] (0,1) -- (1,1);	
	\draw[thick] (1,-1) -- (1,1);	
\end{tikzpicture}\qquad
\begin{tikzpicture}[scale=0.4,x=1cm ]
	\draw[thin,-latex] (-1,-1) -- (3,-1);
	\draw[thin,-latex] (-1,-1) -- (-1,2);		
	
	\draw[thick] (-1,-1) -- (-1,1);
	\draw[thick] (-1,1) -- (0,1);		
	\draw[thick] (0,-1) -- (0,1);		
	\draw[thick] (0,0) -- (1,0);	
	\draw[thick] (1,0) -- (1,-1);	
\end{tikzpicture}
\qquad
\begin{tikzpicture}[scale=0.4,x=1cm ]
	\draw[thin,-latex] (-1,-1) -- (3,-1);
	\draw[thin,-latex] (-1,-1) -- (-1,2);		
	
	\draw[thick] (-1,-1) -- (-1,0);
	\draw[thick] (-1,0) -- (2,0);	
	\draw[thick] (1,0) -- (1,-1);	
		\draw[thick] (2,0) -- (2,-1);	
			\draw[thick] (0,0) -- (0,-1);	
\end{tikzpicture}
	\end{center}
	\caption{Small examples of the bijection}
\end{figure}
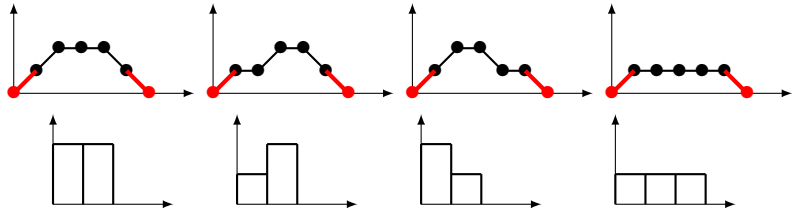


\section{The occurrences of UD and DU are counted.}

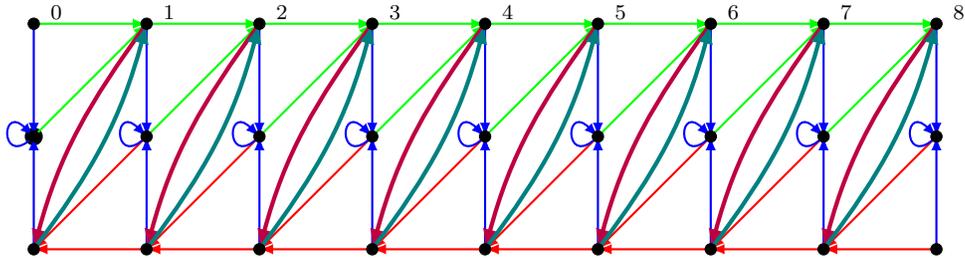
\begin{figure}[h]

	\begin{center}
		\begin{tikzpicture}[scale=1.5,main node/.style={circle,draw,font=\Large\bfseries}]

			\fill (0,-1) circle (0.08cm);

			\foreach \x in {0,1,2,3,4,5,6,7}
			{
				\draw[thick, -latex,green] (\x,0-1) to  (\x+1,0);   
				\draw[ thick, latex-,green] (\x+1,0) to  (\x,0);    
				\draw[ thick,red, -latex] (\x+1,-1) to  (\x,-2);    
				\draw[thick,latex-,blue ] (\x,-1)  ..  controls (\x-0.3,-1+0.3) and  (\x-0.3,-1-0.3) .. (\x,-1) ;	
				\draw[ thick,red, - latex] (\x+1,-2)   to  (\x,-2);    
				
				\node at  (\x+0.2,0.1){\tiny$\x$};
			}  
			
			\foreach \x in {8}
			{
				\draw[thick,latex-,blue ] (\x,-1)  ..  controls (\x-0.3,-1+0.3) and  (\x-0.3,-1-0.3) .. (\x,-1) ;
			}
			
			\foreach \x in {0,1,2,3,4,5,6,7,8}
			{
				
				\draw[ thick,blue, -latex] (\x,0)to  (\x,-1);    
				\draw[ thick,blue, -latex] (\x,-2)to  (\x,-1);    
				
			}                        
			
			\foreach \x in {1,2,3,4,5,6,7,8}
			{
				
				\draw[ ultra thick ,purple, -latex] (\x,0)to  [out=-126,in=75] (\x-1,-2);  
				\draw[ ultra thick ,teal, latex-] (\x,0)to  [out=-105,in=53] (\x-1,-2);      
				
			}

			\node at  (8+0.2,0.1){\tiny$8$};
			
			\foreach \x in {0,1,2,3,4,5,6,7,8}
			{
				\draw (\x,0) circle (0.05cm);
				\fill (\x,0) circle (0.05cm);
				\draw (\x,-1) circle (0.05cm);
				\fill (\x,-1) circle (0.05cm);
				\draw (\x,-2) circle (0.05cm);
				\fill (\x,-2) circle (0.05cm);
			}
		\end{tikzpicture}
	\end{center}
	\caption{Graph (automaton) to count  Motzkin paths with `forbidden' subwords UD and DU.}
	\label{purpelsigma}
\end{figure}

We start in Figure~\ref{purpelsigma} with the relevant automaton. This automaton is based on Figure~\ref{purpel}, but the extra steps are
introduced as ultra thick arcs, which will be labelled by $\sigma z$ resp. $\tau z$.

From the automaton, we deduce sequences of generating functions. Since there are three layers of states, we will
have $f_j, g_j, h_j$ for $j\ge0$. Each of these functions depends on $z,\sigma, \tau$ which is not written, for sake of clarity.
The recursions are obtained by looking at the last step, that leads to the state in question.
\begin{align*}
	f_{i+1}&=zf_i+zg_i+\sigma zh_i,\ f_0=0,\\*
	g_i&=[i=0]+zf_i+zg_i+zh_i,\\*
	h_i&=\tau zf_{i+1}+zg_{i+1}+zh_{i+1}.
\end{align*}
Then we have generating functions
\begin{equation*}
F(u)=\sum_{j\ge0}u^jf_j,\quad G(u)=\sum_{j\ge0}u^jg_j,\quad H(u)=\sum_{j\ge0}u^jh_j,
\end{equation*}
and, by summing the above recursions,
\begin{align*}
	F(u)&=zuF(u)+zuG(u)+\sigma zuH(u),\\*
	G(u)&=1+zF(u)+zG(u)+H(u),\\*
	H(u)&=\frac {\tau z}uF(u)+\frac zu(G(u)-G(0))+\frac zu(H(u)-H(0)).
\end{align*}
To write out the solutions, we use abbreviations:
\footnotesize
\begin{equation*}
W={\sqrt{(1-3 z+ ( 1-\tau \sigma ) {z}^{2}+ ( \tau \sigma-
	\sigma+1-\tau ) {z}^{3})(1+z+ ( 1-\tau \sigma ) {z}^{2}+ ( \tau \sigma-\sigma
	+1-\tau ) {z}^{3})}},
\end{equation*}
\normalsize
and can factorize the (common) denominators; we have $=z(u-r_1)(u-r_2)$ with
\begin{align*}
	r_1&=\frac{1-z+ ( 1-\tau \sigma ) {z}^{2}+ ( \tau \sigma-\sigma
		+1-\tau ) {z}^{3}+W}{2z},\\*
	r_2&=\frac{1-z+ ( 1-\tau \sigma ) {z}^{2}+ ( \tau \sigma-\sigma
		+1-\tau ) {z}^{3}-W}{2z},
\end{align*}
and thus the solutions
\footnotesize
\begin{align*}
F(u)&=\frac{-zu ( \sigma z-z-z\sigma G(0)-z\sigma H(0)+{z}^{2}
	\sigma G(0)+{z}^{2}\sigma H(0)-{z}^{2}G(0)-{z}^{2}H(0)+u )}{z(u-r_1)(u-r_2)},\\
	G(u)&=\frac{z{u}^{2}+{z}^{3}\sigma uG(0)-{z}^{3}uG(0)+{z}^{3}\sigma uH(0)-{z}^{3}uH(0)+\tau \sigma {z}^{2}u-{z}^{2}u-u+{z}^{2}
		G(0)+{z}^{2}H(0)+z
	}{z(u-r_1)(u-r_2)},\\
H(u)&=\frac{-z ( -zu+1+zuH(0)+zu\tau-G(0)-H(0)+zG(0)+zuG(0)+zH(0) )
}{z(u-r_1)(u-r_2)}.
\end{align*}
\normalsize The factor $u-r_2$ can be divided out (this is the `bad' factor, as in many examples of \cite{garden}):\footnote{$1/(u-r_2)$ has no power series expansion around
$u=z=0$, which cannot happen.}
\begin{align*}
	F(u)&=\frac{- ( r_2+u+{z}^{2}\sigma G(0)+{z}^{2}\sigma H(0)-{z}^{2}G(0)-{z}^{2}H(0)+\sigma z-z-z\sigma G(0)-z\sigma H(0) ) z	}{z(u-r_1)},\\
	G(u)&=\frac{r_2z+{z}^{3}\sigma G(0)-{z}^{3}G(0)+{z}^{3}\sigma H(0)-{z}
		^{3}H(0)+\tau \sigma {z}^{2}-{z}^{2}-1+zu
		}{z(u-r_1)},\\
	H(u)&=\frac{-{z}^{2} ( H(0)+\tau-1+G(0) )	}{z(u-r_1)}.
\end{align*}
Now one can set $u=0$ and compute $G(0)$ and $H(0)$; $F(0)=0$ by inspection. Then
\begin{align*}
	G(0)&=\frac{r_2z+{z}^{3}\sigma G(0)-{z}^{3}G(0)+{z}^{3}\sigma H(0)-{z}
		^{3}H(0)+\tau \sigma {z}^{2}-{z}^{2}-1	}{z(-r_1)},\\
	H(0)&=\frac{-{z}^{2} ( H(0)+\tau-1+G(0) )	}{z(-r_1)}=\frac{{z} ( H(0)+\tau-1+G(0) )	}{r_1},
\end{align*}
and further
\begin{align*}
G(0)&={\frac {\tau \sigma {z}^{2}-\tau {z}^{3}\sigma+\tau {z}^{3}+
		\sigma {z}^{3}-{z}^{3}-{z}^{2}+z-2 \sigma z+W-1}{2z ( -\sigma+
		\sigma z-z ) }},\\
H(0)&=	{\frac {-\tau \sigma {z}^{2}+2 \tau {z}^{3}\sigma-\tau {z}^
			{3}-\tau {z}^{4}\sigma+\tau {z}^{4}+\sigma {z}^{4}-{z}^{4}-\sigma 
			{z}^{3}+zW-2 z+1-W}{-2{z}^{2} ( -\sigma+\sigma z-z ) }},\\
	G(0)+H(0)&= {\frac {-1+W-{z}^{3}+\tau{z}^{2}\sigma-\tau{z}^{3}\sigma+\tau
		{z}^{3}+z+{z}^{2}-2 {z}^{2}\sigma+{z}^{3}\sigma}{2{z}^{2}  ( -z-		\sigma+z\sigma  ) }}.
\end{align*}
The expression for $F(u)+G(u)+H(u)$ is somewhat long; we first print  a few terms of the series expansion:
\scriptsize
\begin{gather*}
F(u)+G(u)+H(u)=1+ ( 1+u ) z+ ( {u}^{2}+2 u+\tau+1 ) {z}^{2}+
( 2+\tau \sigma u+{u}^{3}+u\tau+3 u+2 \tau+3 {u}^{2}
) {z}^{3}\\+ ( 4+2 \tau {u}^{2}\sigma+2 \tau \sigma u+4
 {u}^{3}+{u}^{4}+\tau {u}^{2}+4 u\tau+u\sigma+{\tau}^{2}\sigma+5 u
+4 \tau+6 {u}^{2} ) {z}^{4}\\+ ( 8+6 \tau {u}^{2}\sigma+4
 \tau \sigma u+2 u{\tau}^{2}\sigma+u{\tau}^{2}{\sigma}^{2}+10 {u}
^{3}+5 {u}^{4}+6 \tau {u}^{2}+2 {u}^{2}\sigma+\tau {u}^{3}+10 u
\tau\\+2 \tau \sigma+3 u\sigma+{u}^{5}+2 {\tau}^{2}\sigma+3 \tau 
\sigma {u}^{3}+10 u+8 \tau+11 {u}^{2}+{\tau}^{2} ) {z}^{5}\\+{}
( 16+14 \tau {u}^{2}\sigma+2 u{\sigma}^{2}\tau+12 \tau 
\sigma u+8 u{\tau}^{2}\sigma+2 u{\tau}^{2}{\sigma}^{2}+21 {u}^{3}+
\sigma+15 {u}^{4}+19 \tau {u}^{2}\\+8 {u}^{2}\sigma+8 \tau {u}^{3}
+3 {u}^{2}{\tau}^{2}{\sigma}^{2}+2 u{\tau}^{2}+22 u\tau+6 \tau 
\sigma+\tau {u}^{4}+3 {u}^{3}\sigma+7 u\sigma+{\tau}^{3}{\sigma}^{2
}+6 {u}^{5}\\{} +{u}^{6}+6 {\tau}^{2}\sigma+12 \tau \sigma {u}^{3}+4 
\tau \sigma {u}^{4}+3 {\tau}^{2}\sigma {u}^{2}+21 u+18 \tau+22 
{u}^{2}+3 {\tau}^{2} ) {z}^{6}\\ {}+ ( 33+36 \tau {u}^{2}
\sigma+6 u{\sigma}^{2}\tau+36 \tau \sigma u+23 u{\tau}^{2}\sigma+
6 u{\tau}^{2}{\sigma}^{2}+6 \tau {u}^{2}{\sigma}^{2}+44 {u}^{3}+4
 \sigma+36 {u}^{4}\\+2 {\tau}^{3}\sigma+48 \tau {u}^{2}+22 {u}^{2}
\sigma+u{\sigma}^{2}+4 {u}^{4}\sigma+31 \tau {u}^{3}+9 {u}^{2}{	\tau}^{2}{\sigma}^{2}+9 u{\tau}^{2}+50 u\tau+18 \tau \sigma+10 
\tau {u}^{4}\\ {} +3 {\tau}^{2}{u}^{2}+{u}^{5}\tau+15 {u}^{3}\sigma+17 u
\sigma+2 {\tau}^{3}{\sigma}^{2}+21 {u}^{5}+7 {u}^{6}+{u}^{7}+3 u{
	\tau}^{3}{\sigma}^{2}+4 {\tau}^{2}\sigma {u}^{3}+16 {\tau}^{2}
\sigma+3 {\tau}^{2}{\sigma}^{2}\\+33 \tau \sigma {u}^{3}+20 \tau 
\sigma {u}^{4}+{\tau}^{3}{\sigma}^{3}u+18 {\tau}^{2}\sigma {u}^{2}+
5 \tau \sigma {u}^{5}+6 {\tau}^{2}{\sigma}^{2}{u}^{3}+44 u+40 
\tau+47 {u}^{2}+9 {\tau}^{2} ) {z}^{7}+\cdots
	\end{gather*}
\normalsize
Some special cases: Return to the $x$-axis ($u=0$):
	\begin{gather*}
1+z+ ( 1+\tau ) {z}^{2}+2  ( 1+\tau ) {z}^{3}+
( 4+{\tau}^{2}\sigma+4 \tau ) {z}^{4}+ ( 8+2 \tau 
\sigma+2 {\tau}^{2}\sigma+8 \tau+{\tau}^{2} ) {z}^{5}\\+ ( 
16+\sigma+6 \tau \sigma+{\tau}^{3}{\sigma}^{2}+6 {\tau}^{2}\sigma+
18 \tau+3 {\tau}^{2} ) {z}^{6}\\+ ( 33+4 \sigma+2 {\tau}^
{3}\sigma+18 \tau \sigma+2 {\tau}^{3}{\sigma}^{2}+16 {\tau}^{2}
\sigma+3 {\tau}^{2}{\sigma}^{2}+40 \tau+9 {\tau}^{2} ) {z}^{7}+\cdots
	\end{gather*}
$\sigma=0$, $\tau=1$ (peakless paths, final level as coefficient of $u^j$):
\begin{gather*}
1+ ( 1+u ) z+ ( 2 u+{u}^{2}+2 ) {z}^{2}+
( 4+4 u+3 {u}^{2}+{u}^{3} ) {z}^{3}+ ( 8+9 u+7 {u
}^{2}+4 {u}^{3}+{u}^{4} ) {z}^{4}\\ {} + ( 17+20 u+17 {u}^{2}+
11 {u}^{3}+5 {u}^{4}+{u}^{5} ) {z}^{5}+ ( 37+29 {u}^{3}+
45 u+16 {u}^{4}+6 {u}^{5}+41 {u}^{2}+{u}^{6} ) {z}^{6},
\end{gather*}
$\sigma=1$, $\tau=0$ (valleyless)
\begin{gather*}
1+ ( 1+u ) z+ ( 2 u+{u}^{2}+2 ) {z}^{2}+
( 4+4 u+3 {u}^{2}+{u}^{3} ) {z}^{3}+ ( 8+9 u+7 {u
}^{2}+4 {u}^{3}+{u}^{4} ) {z}^{4}\\+ ( 17+20 u+17 {u}^{2}+
11 {u}^{3}+5 {u}^{4}+{u}^{5} ) {z}^{5}+ ( 37+29 {u}^{3}+
45 u+16 {u}^{4}+6 {u}^{5}+41 {u}^{2}+{u}^{6} ) {z}^{6},
\end{gather*}
$\sigma=0$, $\tau=0$ (borderless)
\begin{gather*}
1+ ( 1+u ) z+ ( 1+u ) ^{2}{z}^{2}+ ( 2+u
)  ( {u}^{2}+u+1 ) {z}^{3}+ ( 4+5 u+6 {u}^{2}
+4 {u}^{3}+{u}^{4} ) {z}^{4}\\+ ( 2+u )  ( {u}^{4
}+3 {u}^{3}+4 {u}^{2}+3 u+4 ) {z}^{5}+ ( 16+21 u+22 {u
}^{2}+21 {u}^{3}+15 {u}^{4}+6 {u}^{5}+{u}^{6} ) {z}^{6}\\+
( 33+44 {u}^{3}+36 {u}^{4}+21 {u}^{5}+7 {u}^{6}+{u}^{7}+44 
u+47 {u}^{2} ) {z}^{7}.
\end{gather*}
Return to the $x$-axis, peakless ($u=0$, $\sigma=0$, $\tau=1$) leads to sequence A004148 in \cite{OEIS}
\begin{equation*}
{\frac {{z}^{2}+z-1+\sqrt { ( 1+{z}^{2}+z )  ( {
				z}^{2}-3 z+1 ) }}{-2{z}^{3}}}=
			1+z+2 z^2+4 z^3+8 z^4+17 z^5+37 z^6+82 z^7+\cdots.
\end{equation*}
Return to the $x$-axis, valleyless ($u=0$, $\sigma=1$, $\tau=0$) leads to sequence A004148
\begin{equation*}
{\frac {z-1-{z}^{2}+\sqrt { ( 1+{z}^{2}+z )  ( {
				z}^{2}-3 z+1 ) }}{-2{z}^{2}}}=1+z+z^2+2z^3+4z^4+8z^5+17z^6+37z^7+\cdots.
			\end{equation*}
Return to the $x$-axis, borderless ($u=0$, $\sigma=0$, $\tau=0$) leads to sequence A004149
\begin{equation*}
{\frac {z+\sqrt {(1-z^4) ( -{z}^{2}-2 z+1 ) }-1+{z}^{2}-{
			z}^{3}}{-2{z}^{3}}}=1+z+z^2+2 z^3+4 z^4+8 z^5+16 z^6+33 z^7+\cdots.
\end{equation*}
Ordinary Motzkin paths ($u=0$, $\sigma=1$, $\tau=1$) occur in sequence A001006
\begin{equation*}
{\frac {1-z-\sqrt {1-2z-3z^2}}{2{z}^{2}}}=1+z+2 z^2+4 z^3+9 z^4+21 z^5+51 z^6+127 z^7+323 z^8+\cdots.
\end{equation*}
 For $u=1$, $\sigma=0$, $\tau=1$ (peakless, arbitrary end) we obtain a sequence that is not in \cite{OEIS}
\begin{equation*}
1+2z+5z^2+12z^3+29z^4+71z^5+175z^6+434z^7+1082z^8+2709z^9+6807z^{10}+\cdots.
\end{equation*}
For $u=1$, $\sigma=1$, $\tau=0$ (valleyless, open end) we obtain sequence A091964
\begin{equation*}
1+2z+4z^2+9z^3+21z^4+50z^5+121z^6+296z^7+730z^8+1812z^9+4521z^{10}+\cdots
\end{equation*}
Further, $u=1$, $\sigma=0$, $\tau=0$ (cornerless, open end): sequence A308435
\begin{equation*}
1+2z+4z^2+9z^3+20z^4+45z^5+102z^6+233z^7+535z^8+1234z^9+2857z^{10}+\cdots.
\end{equation*}
Finally we get  open Motzkin paths ($u=1$, $\sigma=1$, $\tau=1$), which is sequence  A005773
 \begin{equation*}
1+2z+5z^2+13z^3+35z^4+96z^5+267z^6+750z^7+2123z^8+6046z^9+17303z^{10}+\cdots.
 \end{equation*}
Here are open Motzkin paths ($u=1$) with $\sigma$ and $\tau$ still unspecified:
\begin{gather*}
F(1)+G(1)+H(1)=1+2 z+  ( 4+\tau  ) {z}^{2}+  ( \tau \sigma+9+3 \tau
 ) {z}^{3}+  ( \sigma+9 \tau+20+4 \tau \sigma+{\tau}^{2}
\sigma  ) {z}^{4}\\* {} +  ( 25 \tau+5 \sigma+4 {\tau}^{2}\sigma+
15 \tau \sigma+45+{\tau}^{2}+{\tau}^{2}{\sigma}^{2}  ) {z}^{5}\\* {} +
 ( 68 \tau+19 \sigma+17 {\tau}^{2}\sigma+5 {\tau}^{2}+48 
\tau \sigma+102+5 {\tau}^{2}{\sigma}^{2}+2 \tau {\sigma}^{2}+{\tau
}^{3}{\sigma}^{2}  ) {z}^{6}\\* {} +  ( 233+148 \tau \sigma+180 
\tau+21 {\tau}^{2}+{\sigma}^{2}+24 {\tau}^{2}{\sigma}^{2}+61 {\tau}
^{2}\sigma+12 \tau {\sigma}^{2}+62 \sigma+{\tau}^{3}{\sigma}^{3}+5
 {\tau}^{3}{\sigma}^{2}+2 {\tau}^{3}\sigma  ) {z}^{7}+\cdots
\end{gather*}
Knowing $F(0), G(0), H(0)$ we can write the full set of solutions (the expressions are somewhat lengthy, and we abbreviate $R=r_2$):
\scriptsize
\begin{equation*}
F(u)=-{\frac {z  ( {R}^{2}z-{z}^{3}R-\tau {z}^{3}\sigma R+\tau {z}^{
			3}R+{z}^{3}\sigma R+R\tau {z}^{2}\sigma-R{z}^{2}+uRz+Rz-R-\tau {z}^
		{3}u\sigma+\tau {z}^{3}u+\tau {z}^{2}u\sigma+zu+z-u  ) }{
		 ( Rz-\tau {z}^{3}\sigma+\tau {z}^{3}+\tau {z}^{2}\sigma-{z}^{
			2}-1+z-{z}^{3}+{z}^{3}\sigma+zu  )   ( -1+\tau {z}^{3}+\tau
		 {z}^{2}\sigma-\tau {z}^{3}\sigma+z+Rz  ) }}
\end{equation*}
\begin{equation*}
	G(u)=\frac{	\genfrac{}{}{0pt}{}{1+uR{z}^{2}-u\tau {z}^{4}\sigma-2 \tau {z}^{2}\sigma+2 
			\tau {z}^{3}\sigma-\tau {z}^{3}+2 \tau {z}^{3}\sigma R-\tau {z}^
			{4}\sigma R+u\tau {z}^{4}-{z}^{3}R+R{z}^{2}-{z}^{4}R+{R}^{2}{z}^{2}}{-
			z+\tau {z}^{3}u\sigma+{z}^{2}u+\tau {z}^{4}R+{z}^{4}\sigma R-{z}^{5
			}\sigma \tau-zu+{z}^{2}-{\tau}^{2}{z}^{5}{\sigma}^{2}-\tau {z}^{4}
			\sigma+{\tau}^{2}{z}^{5}\sigma+{\tau}^{2}{z}^{4}{\sigma}^{2}+\tau {z}
			^{5}{\sigma}^{2}-2 Rz-{z}^{3}\sigma}}{ ( Rz-\tau {z}^{3}\sigma+
		\tau {z}^{3}+\tau {z}^{2}\sigma-{z}^{2}-1+z-{z}^{3}+{z}^{3}\sigma+zu
		)   ( -1+\tau {z}^{3}+\tau {z}^{2}\sigma-\tau {z}^{3}
		\sigma+z+Rz  ) }
\end{equation*}
\begin{equation*}
H(u)=-{\frac {{z}^{2}  ( -\tau+{z}^{3}{\tau}^{2}+{\tau}^{2}{z}^{2}
		\sigma-{\tau}^{2}{z}^{3}\sigma-\tau {z}^{2}+\tau z+R\tau z-\tau {z
		}^{3}+\tau {z}^{3}\sigma-z  ) }{  ( Rz-\tau {z}^{3}\sigma+
		\tau {z}^{3}+\tau {z}^{2}\sigma-{z}^{2}-1+z-{z}^{3}+{z}^{3}\sigma+zu
		 )   ( -1+\tau {z}^{3}+\tau {z}^{2}\sigma-\tau {z}^{3}
		\sigma+z+Rz  ) }}
\end{equation*}
\normalsize
The coefficients of $u^j$ could be written out, as $u$ appears only in linear form in the denominator. 
The coefficients of $u^0$ (return to the origin), viz. $G(0)$ and $H(0)$ have been given just above.

 \section{ Enumeration of skew Motzkin paths relative to peaks and valleys }

As already indicated, we move to the family of skew Motzkin paths \cite{Qing}, \cite{garden}.
There is an additional left-step $L=(-1,-1)$, and $UL$ and $LU$ are forbidden, to avoid overlaps.
The counts of $UD$ and $DU$ are performed with extra variables $\sigma$, $\tau$.
To model this instance, we need a fourth layer, for which we use letters $k_j$ and $K(u)$.
For left-steps we use curvy arcs.
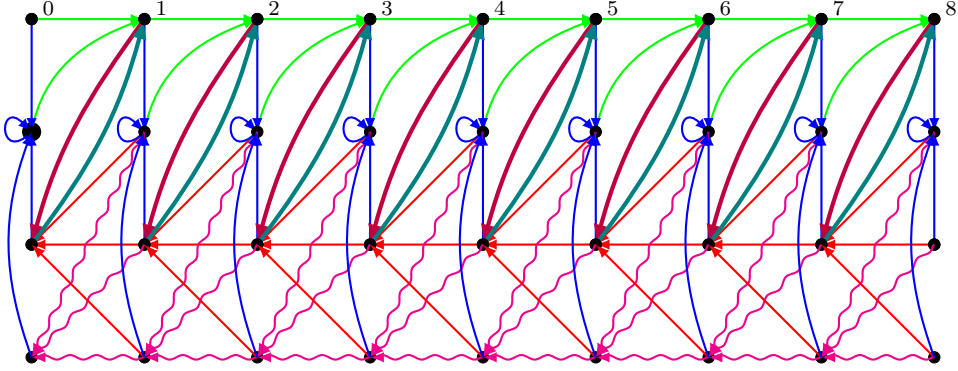
\begin{figure}[h]

	\begin{center}
		\begin{tikzpicture}[scale=1.5,main node/.style={circle,draw,font=\Large\bfseries}]

			\fill (0,-1) circle (0.085cm);

			\foreach \x in {0,1,2,3,4,5,6,7}
			{
				\draw[thick, -latex,green] (\x,0-1) to  [out=80,in=-160](\x+1,0);   
				\draw[ thick, latex-,green] (\x+1,0) to  (\x,0);    
				\draw[ thick,red, -latex] (\x+1,-1) to  (\x,-2);    
				\draw[thick,latex-,blue ] (\x,-1)  ..  controls (\x-0.3,-1+0.3) and  (\x-0.3,-1-0.3) .. (\x,-1) ;	
				\draw[ thick,red, - latex] (\x+1,-2)   to  (\x,-2);    
				
				\node at  (\x+0.15,0.1){\tiny$\x$};
			}  
			
			\foreach \x in {8}
			{
				\draw[thick,latex-,blue ] (\x,-1)  ..  controls (\x-0.3,-1+0.3) and  (\x-0.3,-1-0.3) .. (\x,-1) ;
			}
			
			\foreach \x in {0,1,2,3,4,5,6,7,8}
			{
				
				\draw[ thick,blue, -latex] (\x,0)to  (\x,-1);    
				\draw[ thick,blue, -latex] (\x,-2)to  (\x,-1);    
				
			}                        
			
			\foreach \x in {1,2,3,4,5,6,7,8}
			{
				
				\draw[ ultra thick ,purple, -latex] (\x,0)to  [out=-126,in=75] (\x-1,-2);  
				\draw[ ultra thick ,teal, latex-] (\x,0)to  [out=-105,in=53] (\x-1,-2);      
				
			}

			\node at  (8+0.15,0.1){\tiny$8$};
			
			\foreach \x in {0,1,2,3,4,5,6,7,8}
			{
				\draw (\x,0) circle (0.05cm);
				\fill (\x,0) circle (0.05cm);
				\draw (\x,-1) circle (0.05cm);
				\fill (\x,-1) circle (0.05cm);
				\draw (\x,-2) circle (0.05cm);
				\fill (\x,-2) circle (0.05cm);
			}
		
		\foreach \x in {0,1,2,3,4,5,6,7,8}
		{
			\draw (\x,0) circle (0.05cm);
			\fill (\x,0) circle (0.05cm);
			\draw (\x,-2) circle (0.05cm);
			\fill (\x,-1) circle (0.05cm);
			\draw (\x,-3) circle (0.05cm);
			\fill (\x,-3) circle (0.05cm);
		}
		\foreach \x in {1,2,3,4,5,6,7,8}
	{
		\draw[snake=snake, segment amplitude=.3mm,magenta, thick, -latex] (\x,-2) to (\x-1,-3);
				\draw[snake=snake, segment amplitude=.3mm,magenta, thick, -latex] (\x,-1) to (\x-1,-3);
				\draw[red, thick, -latex] (\x,-3) to (\x-1,-2);

				\draw[snake=snake, segment amplitude=.3mm,magenta, thick, -latex] (\x,-3) to (\x-1,-3);
	}
		\foreach \x in {0,1,2,3,4,5,6,7,8}
{
	\draw[blue, thick, -latex] (\x,-3) to [out=110,in=-110](\x,-1);

}

		\end{tikzpicture}
	\end{center}
	\caption{Graph (automaton) to recognize  skew Motzkin paths with forbidden subwords UD and DU.}
	\label{purpelUDDU2}
\end{figure}

The system of recursions is as follows:
\begin{align*}
	f_{i+1}&=zf_i+zg_i+\sigma zh_i,\ f_0=0,\\*
	g_i&=[i=0]+zf_i+zg_i+zh_i+zk_i,\\*
	h_i&=\tau zf_{i+1}+zg_{i+1}+zh_{i+1}+zk_{i+1},\\
	k_i&=zg_{i+1}+zh_{i+1}+zk_{i+1},
\end{align*}
and as a consequence
\begin{align*}
	F(u)&=zuF(u)+zuG(u)+\sigma zuH(u),\\*
	G(u)&=1+zF(u)+zG(u)+zH(u)+zK(u),\\*
	H(u)&=\frac {\tau z}uF(u)+\frac zu(G(u)-G(0))+\frac zu(H(u)-H(0))+\frac zu(K(u)-K(0)),\\
	K(u)&=\frac zu(G(u)-G(0))+\frac zu(H(u)-H(0))+\frac zu(K(u)-K(0)).
\end{align*}

The system can be solved:
\begin{align*}
F(u)&=\frac{\mathcal{F}}{z(u-r_1)(u-r_2)},\\
	G(u)&=\frac{\mathcal{G}		}{z(u-r_1)(u-r_2)},\\
	H(u)&=\frac{\mathcal{H}}{z(u-r_1)(u-r_2)},\\
	K(u)&=\frac{\mathcal{K}}{z(u-r_1)(u-r_2)};
\end{align*}
the following abbreviations are used:
\begin{gather*}
	\mathcal{F}=-zu ( -2 {z}^{2}K(0)+{z}^{2}\sigma K(0)+{z}^{2}\sigma 
G(0)+{z}^{2}\sigma H(0)-K(0) z\sigma-G(0) z\sigma\\-H(0) z\sigma-2 {z}^{2}G(0)-2 {z}^{2}H(0)+z\sigma-2 z+u
)\\
\mathcal{G}=-{z}^{3}\sigma \tau-{z}^{4}\sigma \tau G(0)-{z}^{4}\sigma \tau
H(0)-{z}^{4}\sigma \tau K(0)+u{z}^{2}\sigma \tau+z{u}^{2}
-2 {z}^{2}u\\+K(0) {z}^{3}u\sigma+G(0) {z}^{3}u\sigma+H(0) {z}^{3}u\sigma-2 {z}^{3}uK(0)-2 {z}^{3}uG(0)\\-2 {z}^{3}uH(0)+2 {z}^{2}G(0)+2 {z}^{2}H(0)+2 {z}^{2}K(0)+2 z-		u\\
\mathcal{H}=( \tau {z}^{2}+{z}^{3}\tau G(0)+{z}^{3}\tau H(0)+{z}^
{3}\tau K(0)-\tau zu+H(0)+G(0)-zK(0)\\+zu-zG(0)-zH(0)-uzG(0)-uzH(0)-uzK(0)-1+K(0) ) z		\\
	\mathcal{K}=z ( {z}^{2}\sigma \tau-{z}^{2}\sigma \tau H(0)-{z}^{2}
	\sigma \tau K(0)+{z}^{3}\sigma \tau G(0)-{z}^{2}\sigma 
	\tau G(0)+H(0)\\+{z}^{3}\sigma \tau H(0)+{z}^{3}\sigma 
	\tau K(0)+G(0)-zK(0)+zu-zG(0)-zH(0)\\-uzG(0)-uz
	H(0)-\tau {z}^{2}-{z}^{3}\tau G(0)-{z}^{3}\tau H(0)\\-uzK(0)-1+K(0)-{z}^{3}\tau K(0) ).		
\end{gather*}
Again, the factors $u-r_2$ can be cancelled out from numerators resp.\ denominators; the square root and
the solutions of the quadratic equation (= the denominator $(u-r_1)(u-r_2)$) are
\begin{gather*}
	W=\bigl(1-2 {z}^{2}\sigma \tau+4 {z}^{3}\sigma \tau-2 {z}^{4}
	\sigma \tau-3 {z}^{2}+2 {z}^{4}\tau+{z}^{6}{\tau}^{2}-4 {z}^{6}
	\tau-4 {z}^{5}\tau+{z}^{6}{\sigma}^{2}-4 {z}^{6}\sigma\\-4 {z}^{5}
	\sigma+4 {z}^{6}+8 {z}^{5}+6 {z}^{6}\tau \sigma-2 {z}^{6}{\tau}^{
		2}\sigma+2 {z}^{5}{\tau}^{2}\sigma-2 {z}^{6}{\sigma}^{2}\tau+2 {z}^
	{5}{\sigma}^{2}\tau\\+{z}^{6}{\sigma}^{2}{\tau}^{2}-2 {z}^{5}{\sigma}^{
		2}{\tau}^{2}+{z}^{4}{\sigma}^{2}{\tau}^{2}-2 {z}^{3}\tau+2 {z}^{4}
	\sigma-2 {z}^{3}\sigma-2 z\bigr)^{1/2},
\end{gather*}
\begin{align*}
	r_1&=\frac{{z}^{3}\sigma \tau-{z}^{2}\sigma \tau-{z}^{3}\sigma-{z}^{3}\tau-z+1+
		2 {z}^{2}+2 {z}^{3}+W}{2z},\\
	r_2&=	\frac{{z}^{3}\sigma \tau-{z}^{2}\sigma \tau-{z}^{3}\sigma-{z}^{3}\tau-z+1+
		2 {z}^{2}+2 {z}^{3}-W
	}{2z}.
\end{align*}
Dividing out $u-r_2$,
\footnotesize
\begin{equation*}
	F(u)=\frac{- ( r_2-2{z}^{2}K(0)+{z}(z-1)\sigma (G(0)+H(0)+K(0)                   )-2{z}^{2}
		G(0)-2{z}^{2}H(0)+z\sigma-2z+u) z
	}{z(u-r_1)},
\end{equation*}
\normalsize
\begin{equation*}
	G(u)=\frac{r_2z+{z}^{3}(\sigma-2)(G(0)+H(0) +K(0)                 )+{z}^{2}\sigma\tau-2{z}^{2}
		-1+zu}{z(u-r_1)},
\end{equation*}
\begin{equation*}
	H(u)=\frac{z^2(-1+\tau+G(0)+H(0)+K(0))}{z(u-r_1)},
\end{equation*}
\begin{equation*}
K(u)=\frac{-z^2(-1+G(0)+H(0)+K(0))}{z(u-r_1)}.
\end{equation*}
The total generating function is
\begin{equation*}
F(u)+G(u)+H(u)+K(u)=\frac{\mathcal A}{\mathcal B},
\end{equation*}
with abbreviations
\begin{gather*}
	\mathcal A=r_2{z}^{3}\sigma \tau-r_2{z}^{3}\tau-r_2z-{z}^{5}{\sigma}^{2}{\tau}
	^{2}+2 {z}^{5}\tau-{z}^{5}{\tau}^{2}-3 {z}^{5}\sigma \tau+2 {z}^{5
	}{\tau}^{2}\sigma+{z}^{5}{\sigma}^{2}\tau-{z}^{4}\sigma \tau\\*+{z}^{4}{
		\sigma}^{2}{\tau}^{2}-{z}^{4}\sigma {\tau}^{2}+2 {z}^{4}\tau-2 {z}^
	{3}\tau+2 {z}^{3}+2 {z}^{3}\sigma \tau-{z}^{3}\sigma-2 {z}^{2}
	\sigma \tau+\tau {z}^{2}-z+1,\\*
	\mathcal B=( r_2z+{z}^{3}\tau+{z}^{3}\sigma-2
	 {z}^{3}-{z}^{3}\sigma \tau+{z}^{2}\sigma \tau-2 {z}^{2}+z-1+zu
	) \\ \quad \quad\quad\quad\times( r_2z+{z}^{3}\tau-{z}^{3}\sigma \tau+{z}^{2}\sigma 
	\tau+z-1 ).
\end{gather*}
The generating function comprising skew Motzkins returning to the $x$-axis is
\footnotesize
\begin{gather*}
F(0)+G(0)+H(0)+K(0)=
1+z+  ( \tau+1  ) {z}^{2}+  ( 3+2 \tau  ) {z}^{3}+
 ( {\tau}^{2}\sigma+7+5 \tau  ) {z}^{4}\\ {}+  ( {\tau}^{2}+
3 \sigma \tau+2 {\tau}^{2}\sigma+17+12 \tau  ) {z}^{5}+
 ( 2 \sigma+3 {\tau}^{2}+9 \sigma \tau+8 {\tau}^{2}\sigma+{
	\sigma}^{2}{\tau}^{3}+41+33 \tau  ) {z}^{6}\\ {}+  ( 8 \sigma+12
 {\tau}^{2}+34 \sigma \tau+24 {\tau}^{2}\sigma+2 {\sigma}^{2}{
	\tau}^{3}+2 \sigma {\tau}^{3}+4 {\sigma}^{2}{\tau}^{2}+103+86 \tau
 ) {z}^{7}\\ {} +  ( 32 \sigma+40 {\tau}^{2}+110 \sigma \tau+
83 {\tau}^{2}\sigma+12 {\sigma}^{2}{\tau}^{3}+6 \sigma {\tau}^{3}+
12 {\sigma}^{2}{\tau}^{2}+5 {\sigma}^{2}\tau+{\tau}^{3}+{\sigma}^{3}
{\tau}^{4}+259+233 \tau  ) {z}^{8}+   \cdots.
\end{gather*}
\normalsize
The generating function comprising open skew Motzkins (ending at any level) is
\footnotesize
\begin{gather*}F(1)+G(1)+H(1)+K(1)=
1+2 z+  ( 4+\tau  ) {z}^{2}+  ( 3 \tau+\tau \sigma+10
 ) {z}^{3}\\ {} +  ( 24+4 \tau \sigma+10 \tau+\sigma {\tau}^{2
}+\sigma  ) {z}^{4}+  ( 16 \tau \sigma+5 \sigma+30 \tau+
60+4 \sigma {\tau}^{2}+{\tau}^{2}+{\sigma}^{2}{\tau}^{2}  ) {z}
^{5}\\ {}+  ( 5 {\sigma}^{2}{\tau}^{2}+53 \tau \sigma+152+21 \sigma
+5 {\tau}^{2}+90 \tau+19 \sigma {\tau}^{2}+2 {\sigma}^{2}\tau+{
	\tau}^{3}{\sigma}^{2}  ) {z}^{6}\\ {}+  ( 392+{\tau}^{3}{\sigma}^{
	3}+5 {\tau}^{3}{\sigma}^{2}+2 {\tau}^{3}\sigma+24 {\tau}^{2}+12 {
	\sigma}^{2}\tau+75 \sigma+{\sigma}^{2}+72 \sigma {\tau}^{2}+25 {
	\sigma}^{2}{\tau}^{2}+178 \tau \sigma+262 \tau  ) {z}^{7}+ \cdots.
\end{gather*}
\normalsize
For  $\sigma=\tau=1$, we obtain sequence A082582 in \cite{OEIS} and more special cases could be listed since the 
functions $F(u), G(u), H(u), K(u)$ are explicitly known.

\section{Conclusion}

Asymptotic questions are not addressed in this paper. For those interested in the subject,
it is mentioned that everything is driven by square-root singularities; what we abbreviated by $W$, 
leads to a dominated singularity $\rho$. For questions like the average endpoint, a differentiation w.r.t. $u$, 
followed by $u=1$, will be crucial.


\bibliographystyle{plain}


\end{document}